\documentclass[a4paper, 12pt]{amsart}
\usepackage[T1]{fontenc}
\usepackage{amssymb}
\usepackage[utf8]{inputenc}
\usepackage[pagewise]{lineno}
\usepackage{multicol}
\usepackage{hyperref}
\usepackage{bm}
\usepackage{comment}
\usepackage{xcolor}

\makeatletter
\@namedef{subjclassname@2020}{%
  \textup{2020} Mathematics Subject Classification}
\makeatother

\newtheorem{theorem}{Theorem}[section]
\newtheorem{proposition}[theorem]{Proposition}
\newtheorem{lemma}[theorem]{Lemma}
\newtheorem{corollary}[theorem]{Corollary}
\theoremstyle{definition}
\newtheorem{remark}[theorem]{Remark}

\numberwithin{equation}{section}

\begin{document}

\title{Graded identities of the adjoint representation of $\mathfrak{sl}_2(K)$}

\author[Sampaio and Koshlukov]{Cássia Sampaio and Plamen Koshlukov}
\address{Department of Mathematics, State University of Campinas, 651 Sérgio Buarque de Holanda, 13083-859 Campinas, SP, Brazil}
\email{cfsampaio@ime.unicamp.br, plamen@unicamp.br}

\thanks{C. Sampaio was financed in part by the Coordenação de Aperfeiçoamento de Pessoal de Nível Superior - Brasil (CAPES) - Finance Code 001.\\
\indent P. Koshlukov was partially supported by CNPq Grant 307184/2023-4, and by FAPESP grant 2018/23690-6}



\begin{abstract}
Let $K$ be a field of characteristic zero and let $\mathfrak{sl}_2 (K)$ be the 3-dimensional simple Lie algebra over $K$. In this paper we describe a finite basis for the $\mathbb{Z}_2$-graded identities of the adjoint representation of $\mathfrak{sl}_2 (K)$, or equivalently, the $\mathbb{Z}_2$-graded identities for the pair $(M_3(K), \mathfrak{sl}_2 (K))$. We work with the canonical grading on $\mathfrak{sl}_2 (K)$ and the only nontrivial $\mathbb{Z}_2$-grading of the associative algebra $M_3(K)$ induced by that on $\mathfrak{sl}_2(K)$.

\end{abstract}

\maketitle

\section*{Introduction}

Determining a basis of polynomial identities for the matrix algebra of  order $n$ is a challenging problem that draws the attention of many researchers in the theory of algebras with polynomial identities. So far, there have been some results in this direction, but the problem, in general, remains widely open.  

 A finite basis for the identities of the matrix algebra of order 2, $M_2(K)$, over a field of characteristic zero was described by Razmyslov \cite{razmyslov1973existence} and subsequently minimized by Drensky \cite{drenski1981minimal}. In positive characteristic, it is known that the Lie algebra $M_2(K)$, $\text{char}(K)=2$, does not possess a finite basis of identities when $K$ is infinite, a theorem proved by Vaughan-Lee in \cite{vaughan1970varieties}. If $\text{char}(K) > 2$ and $K$ is infinite, Koshlukov \cite{koshlukov2001basis} exhibited a finite basis of identities for $M_2(K)$. To achieve this, he crucially used the invariant theory of the orthogonal group, see \cite{de1982characteristic}, and weak identities \cite{razmyslov1994identities}. Concerning the algebra $M_n(K)$, $n\geq 3$, little is known about its polynomial identities. For a finite field $K$, Genov and Siderov \cite{genov1981basis, genov1982basis}, determined bases of identities for $M_3(K)$ and $M_4(K)$.  In the case of infinite fields, weak identities and graded identities are expected to be an essential tool, as it was the case for $M_2(K)$, and also in the proof that the variety of associative algebras over fields of characteristic zero satisfies the Specht property, the celebrated theorem of Kemer.

Weak identities first appeared in  Razmyslov's work \cite{razmyslov1973existence}, where, in addition to exhibiting a finite basis of identities for $sl_2(K)$ and $M_2(K)$, $\text{char}(K) = 0$, he proved that the variety $\mathrm{var}(M_2(K))$ satisfies the Specht property. The notions of weak identity and identity of a representation might be easily merged in the following sense. Let  $ K  \langle X \rangle $ be the free associative algebra over K freely generated by a countable set  $X = \{ x_1, x_2, \ldots, x_n, \ldots \}$.
If $L$ is a Lie algebra and $\rho \colon L \rightarrow \mathfrak{gl}(V)$ is a Lie algebra representation of $L$ on a vector space $V$, a \emph{polynomial identity of the representation } $\rho$ is a polynomial $f(x_1, \ldots , x_n) \in K \left< X \right>$ 
 such that  $f(\rho (e_1), \ldots , \rho (e_n)) = 0$ for every $e_1$, \dots $e_n \in L$.
Denoting by $M$ the associative algebra generated by $\rho (L)$, the problem of determining the identities of the representation $\rho$ of $L$ is equivalent to finding the identities for the associative-Lie pair $(M,\rho (L))$.

Identities of this type were also crucial in finding a basis of identities for the Lie algebra $sl_2(K)$ when the field $K$ is infinite and of characteristic greater than 2, \cite{vasilovskii1989basis}. In the same context, Koshlukov  \cite{koshlukov1997weak} determined a basis of identities for the natural representation $\mathfrak{sl}_2(F) \rightarrow M_2(F)$. Here we recall that weak identities can be defined in a similar way for other classes of algebras that are not "very far" from associative ones, and they play an important role in studying polynomial identities for Lie, alternative and Jordan algebras.

We recall the notion of a graded identity. As we work with gradings by $\mathbb{Z}_2$, the group of order 2, we consider gradings by this group only. If $A= A_0 \bigoplus A_1$ is a $\mathbb{Z}_2$-graded algebra, then $f(y_1, \ldots , y_m, z_1, \ldots , z_n)$ is a $\mathbb{Z}_2$-graded identity of $A$ if $f(a_1, \ldots , a_m, b_1 , \ldots , b_n)=0$ for every $a_i \in A_0 $ and $b_i \in A_1$. Here we consider $X=Y\cup Z$, $Y$ and $Z$ disjoint infinite sets, and the $\mathbb{Z}_2$-grading on $K\langle X\rangle$ is given by declaring the variables from $Y$ of degree 0 (even ones), and those from $Z$ of degree 1 (odd ones). 

Graded identities have also proven to be an important tool in studying PI-algebras. They were essential, for example, in establishing that the variety of associative algebras satisfies the Specht property in characteristic 0, \cite{kemer1987finite}. We emphasize that the relation between graded identities and ordinary ones established by Kemer (Theorem 1.1, \cite{kemer1991ideals}) is not completely valid in positive characteristic. We note that if two algebras satisfy the same graded identities they also satisfy the same ordinary identities. But the relationship between graded and ordinary identities remains not that clear: for example, the knowledge of the graded identities of an algebra does not give much information about the ordinary ones satisfied by the same algebra. 

Finite bases of identities for the $\mathbb{Z}_2$-graded algebras  $M_2(K)$, with $\text{char} K = 0$, and for $M_{1,1}(G)$, where $G$ is the Grassmann algebra over a field of characteristic 0 were found by Di Vincenzo \cite{di1992graded}. Later on, these results were extended to infinite fields of characteristic different from 2, see \cite{koshlukov2002graded}.
Once again, over an infinite field $K$, with $\mathrm{char} K \neq 2$,  a basis for the $\mathbb{Z}_2$-graded identities of $\mathfrak{sl}_2(K)$ was described in \cite{koshlukov2008graded}. In order to obtain a finite basis of identities for the pair $(M_2(K), gl_2(K))$, Koshlukov and Krasilnikov \cite{koshlukov2012basis} found a basis for the $\mathbb{Z}_2$-graded identities for the pair $(M_2(D), gl_2(D))$, where $D$ is an infinite integral domain. Under the same assumptions for $D$, C\'odamo and Koshlukov \cite{codamo2024specht} determined a basis of weak identities for the $\mathbb{Z}_2$-graded pair $\left(M_2(D), \mathfrak{sl}_2(D)\right)$ and proved that this pair possesses the Specht property. Trishin, in \cite{trishin1983identities}, found finite bases of identities for the irreducible representations of $\mathfrak{sl}_2(K)$ under the assumption that $K$ is an algebraically closed field of characteristic 0. 

Here we combine ideas from \cite{koshlukov2012basis} and \cite{trishin1983identities} in order to describe a basis for the $\mathbb{Z}_2$-graded identities of the adjoint representation of   $\mathfrak{sl}_2(K)$, the 3-dimensional simple Lie algebra over a field of characteristic zero. 
The image  $\mathrm{ad} (sl_2(K) ) $ of this representation generates $M_3(K)$ as an associative algebra. Hence, equivalently, we describe a basis of the graded identities for the pair $(M_3(K), \mathfrak{sl}_2 (K))$. Clearly, if the field $K$ is not algebraically closed there may exist non-isomorphic 3-dimensional simple Lie algebras over $K$. But a standard argument, using extensions of the scalars, shows that these will satisfy the same polynomial identities. Therefore, without loss of generality, one may assume $K$ is algebraically closed.

\section{Preliminaries}

Throughout this paper, $K$ stands for a field of characteristic zero, and all algebras and vector spaces are considered over $K$. As we stated above, we can (and will) assume $K$ is algebraically closed. We use the additive notation for the cyclic group $\mathbb{Z}_2 = \{0,1\}$ of order two. If $A$ is an associative algebra we denote by $A^{(-)}$ the vector space of $A$ together with the Lie bracket operation $[a,b]=ab-ba$. It is well known that $A^{(-)}$ is a Lie algebra. We shall use the bracket notation for the product in a Lie algebra. The commutators without inner brackets are supposed to be left-normed, that is $[x_1, x_2, \ldots , x_n] = [[x_1, x_2, \ldots x_{n-1}], x_n]$ for every $n\ge 3$.

An algebra $A$ is said to be a $\mathbb{Z}_2$\emph{-graded algebra} (also called simply \emph{graded algebra} or, in the associative case, \emph{superalgebra}) if $A = A_0 \bigoplus A_1$, where $A_0$ and $A_1$ are vector subspaces, respecting $A_i A_j \subseteq A_{i+j}$, $i,j \in \mathbb{Z}_2$. We call $(A,L)$ a $\mathbb{Z}_2$-\emph{graded pair} if $L$ is a Lie subalgebra of $A^{(-)}$ such that $L=L_0 \bigoplus L_1$, $L_i = A_i \cap L$, $i=0$, 1. 
Since we consider gradings by $\mathbb{Z}_2$ only we will omit the group and we will often call these simply graded algebras or graded pairs.

Let $K \left< X \right>$ be the free associative algebra freely generated by 
$X = Y \cup Z$, where $Y = \{ y_1, y_2, \ldots , y_m , \ldots \}$ and $Z = \{ z_1, z_2, \ldots , z_n , \ldots \}$ are countable disjoint sets.
We call even (respectively odd) variables the ones from $Y$ (respectively $Z$). Setting the homogeneous degree of a monomial $m \in K\left<X \right>$ as 0 if its degree with 
respect to the variables $z \in Z$ is even and as 1 otherwise, the algebra $K \left< X \right>$ 
is $\mathbb{Z}_2$-graded with the natural grading $K\left< X \right> = K \left< X \right> _0 \bigoplus K \left< X \right> _1$, where $K \left< X \right> _i$ is the vector subspace spanned by the monomials of homogeneous degree $i$, $i=0$, 1. 
Let $L \left< X \right> $ be the Lie subalgebra of $K\left<X \right> ^{(-)}$ generated by $X$. It is well known (Witt Theorem) that $L \left< X \right>$ is the free Lie algebra freely generated by $X$. With the grading inherited from  $K\left<X \right>$, $ L \left< X \right>$ is the free  $\mathbb{Z}_2$-graded Lie algebra freely generated by $X$. Consequently, the pair $(K\left<X \right> , L \left< X \right> )$ is the free graded pair of countable rank in the class of all graded associative-Lie pairs. 
Hereafter, $K\left<X \right>$ will always denote the free $\mathbb{Z}_2$-graded algebra.

We call an ideal $I$ of $K\left<X \right>$ a $T_2$-\emph{ideal} if $I$ is invariant under all graded endomorphisms of $K\left<X \right>$, that is, if $\phi (I) \subseteq I$ for every endomorphism $\phi$ of $K\left<X \right>$ such that $\phi (y_i) \in K\left<X \right> _0$ and $\phi (z_i) \in K\left<X \right> _1$. An ideal $I$ is called a \emph{weak} $T_2$\emph{-ideal}
if it is invariant under the graded endomorphisms of the pair $(K\left<X \right> , L \left< X \right>)$.
As in the ordinary case, we can define a (weak) $T_2$-ideal generated by a set and 
a basis of (weak) graded identities. 

 A polynomial $f (y_1, \ldots , y_m, z_1 , \ldots z_n)\in K\langle X\rangle$ is said to be a \emph{weak graded identity} of the  $\mathbb{Z}_2$-graded pair $(A,L)$ if $f(a_1^{(0)},\ldots , a_m^{(0)}, $ \linebreak $a_1 ^{(1)}, \ldots , a_n^{(1)}) = 0 $ for all $a_1^{(0)}, \ldots , a_m^{(0)} \in L_0$ and $ a_1^{(1)}, \ldots , a_n ^{(1)} \in L_1$. It can be easily proved that the set $T_2(A, L)$ of weak graded identities of a pair $(A, L)$ is a weak $T_2$-ideal in $K\left<X \right>$. The converse also holds, its proof repeats word by word that of the case of ordinary (associative, ungraded) T-ideals.

Let $S\subset K\langle X\rangle$ be a nonempty set of polynomials, then the ideal of graded weak identities generated by $S$ is the least ideal of graded weak identities $I$ that contains $S$. Recall that in this case $S$ is called a basis of $I$.  A polynomial $f\in K\langle X\rangle$ is a consequence of the identities of $S$ whenever $f\in I$. Accordingly, $f$, $g\in K\langle X\rangle$ are equivalent as weak graded identities when they generate the same ideal of weak graded identities. (In other words, $f$ is a consequence of $g$ and $g$ is a consequence of $f$.)

It is well known that over an infinite field every ideal of identities is generated by its multihomogeneous elements. When the base field is of characteristic 0, the ideals of identities are generated by the multilinear polynomials they contain. The proofs of these facts are standard, moreover they transfer word by word to our case.

\section{
\texorpdfstring{Graded identities for the pair $(M,S)$ }{Graded identities for the pair (M,S)}}

Let $\mathfrak{sl}_2 (K)$ be the simple 3-dimensional Lie algebra over the field $K$ and denote by $S$ the image of the adjoint representation of
 $\mathfrak{sl}_2 (K)$, $ \mathrm{ad} \colon g \in \mathfrak{sl}_2 \mapsto \mathrm{ad} g \in \mathfrak{gl} (\mathfrak{sl}_2 (K))$. Then $\mathfrak{gl}(\mathfrak{sl}_2 (K)) \cong M_3(K)$, since the action of  $\mathfrak{sl}_2$ on the 3-dimensional vector space $\mathfrak{sl}_2(K)$ is irreducible. Hence the image of $\mathfrak{sl}_2(K)$ generates, as an associative algebra, the whole algebra of linear transformations in dimension 3. In this way we can write simply $ \mathrm{ad} \colon \mathfrak{sl}_2 (K) \rightarrow M_3(K)$.

Since $\mathfrak{sl}_2 (K)$ is a simple Lie algebra, the adjoint representation is faithful. Therefore,  $S \cong \mathfrak{sl}_2 (K)$ is also a simple Lie algebra.
Taking the canonical basis of $\mathfrak{sl}_2, \, \{h, x, y \}$, with $[h,x]=2x, \, [h,y]=-2y$ and $[x,y]=h$, we have

\begin{displaymath}
S=\left\{ \left(\begin{array}{ccc}
0 & -n & m \\
-2m & 2a & 0 \\
2n & 0 & -2a
\end{array}\right)  \, | \,a, m, n \in K\right\}.
\end{displaymath}

Note that $S$  generates $M_3(K)$ as an associative algebra. We have the following natural $\mathbb{Z}_2$-grading of the Lie algebra $S$.
\begin{displaymath}
S_0=\left\{ \left(\begin{array}{ccc}
0 & 0 & 0 \\
0 & 2a & 0 \\
0 & 0 & -2a
\end{array}\right)\right\} \, \mathrm{and} \,  S_1=\left\{ \left(\begin{array}{ccc}
0 & -m & n \\
-2n & 0 & 0 \\
2m & 0 & 0
\end{array}\right)\right\},
\end{displaymath}
$a$, $m$, $n \in K$. It is well known that, up to a graded isomorphism, the above is the only nontrivial $\mathbb{Z}_2$-grading on $\mathfrak{sl}_2(K)$, and on $S$. 

Denote $M:=M_3(K)$ and consider the
 $\mathbb{Z}_2$-grading exhibited below, which is, up to a graded isomorphism, the only nontrivial $\mathbb{Z}_2$-grading of $M_3(K)$.
\begin{displaymath}
M_0=\left\{ \left(\begin{array}{ccc}
K & 0 & 0 \\
0 & K & K \\
0 & K & K
\end{array}\right) \right\} \,  \, \,  \, \,\mathrm{and} \,  \, \, \, \, M_1=\left\{ \left(\begin{array}{ccc}
0 & K & K \\
K & 0 & 0 \\
K & 0 & 0
\end{array}\right) \right\}.
\end{displaymath}
We draw the readers' attention that the latter grading is compatible with that on $S$: we have that $S_i = M_i \cap S$. This implies $(M, S)$ is a graded pair. 

Therefore, instead of considering the graded pair $(M_3(K), \mathfrak{sl}_2(K))$ we can study the graded pair $(M,S)$.

We denote by $s_n(x_1,\ldots,x_n)=\sum_{\sigma\in S_n} (-1)^\sigma x_{\sigma(1)}\cdots x_{\sigma(n)}$, the \textsl{standard polynomial} of degree $n$, where $S_n$ is the symmetric group of all permutations of $(1,2,\ldots,n)$, and $(-1)^\sigma$ stands for the sign of $\sigma\in S_n$.

\begin{lemma} \label{lemmaidbase}
    The following polynomials are weak graded identities for the graded pair $(M, S)$.
  \begin{equation}
        [y_1, y_2], \label{idcolchy1y2}
    \end{equation}
    \begin{equation}
        y_1zy_2, \label{idyzy}
    \end{equation}
    \begin{equation}
        z_1 y z_2 + z_2 y z_1 ,\label{idzyz}
    \end{equation}
    \begin{equation}
    s_3(z_1, z_2, z_3), \label{s3}
    \end{equation} 
    \begin{equation}
  y_1 z_1 z_2 y_2 - y_2 z_1 z_2 y_1,  \label{2y2z}
    \end{equation}
    \begin{equation}
        z_1y_1y_2z_2-z_2y_1y_2z_1,\label{z_1y_1y_2z_2-z_2y_1y_2z_2}
    \end{equation}
 \begin{equation}
     s_3 (z_1,z_2, y) x - x s_3 (z_1,z_2, y) , \hspace{0,4cm} x = y_i, z_j. \label{s3zzy}
 \end{equation}
    \begin{equation}
        y[z_1,z_2]z_3+z_3z_1yz_2 .\label{id y[z_1,z_2]z_3+z_3z_1yz_2}
    \end{equation}
\end{lemma}

\begin{proof}
    The identities (\ref{idcolchy1y2}), (\ref{idyzy}), (\ref{idzyz}),  (\ref{2y2z}), and \eqref{z_1y_1y_2z_2-z_2y_1y_2z_2}  are verified by direct and easy computation. For example, two diagonal matrices commute, and hence (\ref{idcolchy1y2}) follows. Substituting distinct elements of $S_0$ and of $S_1$ in (\ref{idyzy}), (\ref{idzyz}),  (\ref{2y2z}), \eqref{z_1y_1y_2z_2-z_2y_1y_2z_2}, and \eqref{id y[z_1,z_2]z_3+z_3z_1yz_2}, one gets that all of them vanish.

To see that  (\ref{s3})  is an identity, note that $S_1$ is a 2-dimensional subspace and $s_3(z_1, z_2, z_3)$ is a skew-symmetric polynomial in three odd variables. Therefore substituting in elements from the basis of $S_1$ we have to repeat some element, and the skew-symmetry yields the polynomial vanishes.  The identity \eqref{s3zzy} follows from the fact that $s_3(z_1,z_2, y)$ is a multiple of the Casimir element of the adjoint representation, which belongs to the centre of the universal enveloping algebra of $\mathfrak{sl}_2 (K)$. 
\end{proof}

As immediate consequences of \eqref{idyzy} and \eqref{idzyz} we have the following  weak $\mathbb{Z}_2$-graded identities of $(M,S)$
\begin{align}
     [y,z_1, z_2]&-yz_1z_2-z_2z_1y, \label{[y,z,z]}\\
[z,y_1,y_2] &- zy_1y_2-y_2y_1z    \label{[z,y,y]}\\
        [y_1, z, y_2] -[y_2,z,y_1], & \quad
    [z_1, y, z_2] - [z_2, y, z_1] \label{[z_1, y, z_2] = [z_2, y, z_1]}
\end{align}

Rewriting the standard polynomial of degree 3 in different ways, we have the following identities which are equivalent to \eqref{s3}:
    \begin{align} \label{s3emcolchetes}
        & z_1[z_2, z_3] -z_2[z_1,z_3]+z_3[z_1,z_2]\\
     \label{s3embolinha}
        &(1- \tau _{13})z_1z_2z_3 - [z_1, z_3] \circ z_2,
    \end{align}
where $\tau _{13}$ is the transposition of the variables $z_1$ and $z_3$, and the symbol $\circ$ stands for the Jordan product, that is $a \circ b = ab+ba$.

Let $I$ be the weak $T_2$-ideal generated by the polynomials of the previous Lemma. We will prove that $I = T_2 (M,S)$.
\begin{lemma}
    The following weak graded identities are valid in the pair $(M,S)$.
\begin{equation}
    z_1yz_2z_3=z_3[z_2,z_1]y \label{z_1yz_2z_3+z_3[z_1,z_2]y}
\end{equation}
    \begin{equation}
        y z_1 z_2z_3= [y, z_2, z_3]z_1 \label{ y z_1 z_2z_3= [y, z_2, z_3]z_1}
    \end{equation}
    \begin{equation}
        z_1 z_2 z_3y = z_3[y, z_1, z_2] \label{ z_1 z_2 z_3y = z_3[y, z_1, z_2]}
    \end{equation}
\end{lemma}
\begin{proof}
To prove \eqref{z_1yz_2z_3+z_3[z_1,z_2]y} we use \eqref{idzyz}, \eqref{s3zzy} and \eqref{id y[z_1,z_2]z_3+z_3z_1yz_2}
\begin{align*}
    0 & \stackrel{\eqref{s3zzy}}{=} z_3(y[z_1,z_2] + z_2yz_1) - (y[z_1,z_2] + z_2yz_1)z_3 \\
    &\stackrel{\eqref{idzyz},\eqref{id y[z_1,z_2]z_3+z_3z_1yz_2}}{=}  z_1yz_2z_3-z_3[z_2,z_1]y.
\end{align*}
With the following manipulations, we obtain \eqref{ y z_1 z_2z_3= [y, z_2, z_3]z_1}.
    \begin{align*}
         [y, z_2, z_3] z_1& \stackrel{\eqref{[y,z,z]}}{=} y z_2 z_3 z_1 +   z_3z_2yz_1 \stackrel{\eqref{id y[z_1,z_2]z_3+z_3z_1yz_2}}{=} y z_2 z_3 z_1 +y[z_1,z_2]z_3 \\
         &= y z_2 [z_3, z_1] +y z_1 z_2 z_3 \stackrel{\eqref{idyzy}}{=}y z_1 z_2 z_3 
    \end{align*}
    Similarly, we prove \eqref{ z_1 z_2 z_3y = z_3[y, z_1, z_2]}.
\end{proof}

\begin{lemma}
The following identity holds in the pair $(M,S)$
    \begin{align}
        &[y, z_{11},z_{12}, z_{21}, z_{22}, \ldots , z_{k1}, z_{k2}] \nonumber\\
        &= [y, z_{\sigma(1)1},z_{\sigma(1)2}, z_{\sigma(2)1}, z_{\sigma(2)2}, \ldots , z_{\sigma(k)1}, z_{\sigma(k)2}], \label{idgeraltrocadoszemcolchete}
    \end{align}
    for every $k>1$ and for every permutation  $\sigma \in S_k$ of the indices 1, \dots, $k$.
\end{lemma}
In particular, we have
\begin{equation}
    [y, z_1, z_2, \ldots , z_{2k-1}, z_{2k}]=[y, z_{2k-1}, z_{2k}, z_{2k-3}, z_{2k-2}, \ldots , z_1 , z_2] .\label{idtrocadoszemcolchete}
\end{equation}
\begin{proof}
    When $k=2$, the identity \eqref{idgeraltrocadoszemcolchete} is equivalent to $[y,z_1, z_2, z_3, z_4] = [y, z_3, z_4, z_1, z_2]$. Now we compute $ [y,z_1, z_2, z_3, z_4]$, it equals
    \begin{align*}
        \stackrel{\eqref{[y,z,z]}, \eqref{idyzy}}{=} & yz_1z_2z_3z_4+z_4z_3yz_1z_2+z_2z_1yz_3z_4+z_4z_3z_2z_1y\\
        & +yz_3[z_4,z_1]z_2+z_2[z_1,z_4]z_3y \\
        \stackrel{\eqref{idyzy}}{=} & yz_3z_4z_1z_2+z_2z_1z_4z_3y +z_4z_3yz_1z_2 -yz_3z_1z_4z_2-z_2z_4z_1z_3y\\
        &+yz_1z_3z_2z_4+z_4z_2z_3z_1y \\
       = &[y, z_3,z_4,z_1,z_2] -y[z_3,z_1]z_4z_2-yz_1z_3z_4z_2 - z_2z_4[z_1,z_3]y \\
        &-z_2z_4z_3z_1y+ yz_1z_3z_2z_4+z_4z_2z_3z_1y \\
         \stackrel{}{=} & [y, z_3,z_4,z_1,z_2]+y[z_1,z_3]z_4z_2+yz_1z_3[z_2,z_4]\\
         &-z_2z_4[z_1,z_3]y-[z_2,z_4]z_3z_1y \\
         \stackrel{\eqref{y2z2=z2y2}, \eqref{2y2z}}{=}&[y, z_3,z_4,z_1,z_2]. 
    \end{align*}
    To prove the general case, we apply iteratively the case proved above, using the fact that $[y,z_{i_1}, z_{i_2}, \ldots, z_{i_{2l}}] \in L \left< X \right> _0$. 
    \end{proof}

\begin{lemma}
    The following polynomials belong to $I$, for every non negative integer $k$,
     \begin{equation}
        y_1 z_1 \ldots z_{2k+1}y_2 \label{zimpar}
    \end{equation}
    \begin{equation}
        z_1 y_1 \ldots y_{2k+1}z_2+z_2y_1 \ldots y_{2k+1}z_1 .\label{yimpar}
    \end{equation}
    \begin{equation}
     z_1y_1y_2 \ldots y_{2k}z_2-z_2y_1y_2 \ldots y_{2k}z_1   \label{zy2z}
    \end{equation} 
     \begin{equation}
        y_1 z_1 \ldots z_{2k}y_2 - y_2 z_1 \ldots z_{2k}y_1, \label{idyz2ky}
    \end{equation} 
\end{lemma}

\begin{proof}

First, we note that by (\ref{idyzy}), the identity \eqref{zimpar}  holds for $k=0$. Since $[ y_2,z_i, z_j] \in S_0$ we have by induction
\begin{align*}
    0= \, & y_1 z_1z_2 \ldots z_{2k-1} [ y_2, z_{2k+1}, z_{2_k}] \\
    \stackrel{\eqref{[y,z,z]}}{=}& y_1 z_1z_2\ldots z_{2k-1}y_2  z_{2k+1}z_{2_k}+ y_1 z_1z_2\ldots z_{2k-1}z_{2k} z_{2k+1}y_2 \\
=\, &y_1 z_1z_2\ldots z_{2k-1}z_{2k} z_{2k+1}y_2 .
\end{align*}
  Similarly, we prove \eqref{yimpar} using  \eqref{idzyz}, the fact that  $[z, y_i, y_j] \in S_1$ and the weak graded identities \eqref{idcolchy1y2}, \eqref{idyzy} and \eqref{[z,y,y]}.

When $k=1$, \eqref{zy2z} is equal to  \eqref{z_1y_1y_2z_2-z_2y_1y_2z_2}. Using \eqref{idcolchy1y2}, \eqref{[z,y,y]}, and induction, we prove \eqref{zy2z} for every $k$. 

Finally, we use \eqref{ y z_1 z_2z_3= [y, z_2, z_3]z_1} and \eqref{idtrocadoszemcolchete} to prove \eqref{idyz2ky}, observing that \eqref{idyz2ky} is equal to \eqref{2y2z}, for $k=1$.
\begin{align*}
    y_1z_1 z_2 \ldots z_{2k} y_2& \stackrel{\eqref{ y z_1 z_2z_3= [y, z_2, z_3]z_1}}{=} [y_1, z_2, z_3]z_1 z_4 z_5 \ldots z_{2k} y_2 \\
    & \vdots \\ 
     & \stackrel{\eqref{ y z_1 z_2z_3= [y, z_2, z_3]z_1}}{=} [y_1, z_2, z_3, \ldots , z_{2k-2}, z_{2k-1}]z_1  z_{2k} y_2 \\
     & \stackrel{\eqref{2y2z}}{=} y_2 z_1 [y_1, z_2, z_3, \ldots , z_{2k-2}, z_{2k-1}] \\
     & \stackrel{\eqref{idtrocadoszemcolchete}}{=}y_2 z_1 [y_1,z_{2k-2}, z_{2k-1} , \ldots ,  z_2, z_3] \\
     & \stackrel{\eqref{ y z_1 z_2z_3= [y, z_2, z_3]z_1}}{=} y_2 z_1 \ldots z_{2k} y_1.\qedhere
\end{align*}
\end{proof}

\begin{lemma}
    The following weak graded identity holds in the graded pair $(M, S)$
    \begin{equation}
        y_1y_2z_1z_2-z_1z_2y_1y_2. \label{y2z2=z2y2}
    \end{equation}
\end{lemma}
\begin{proof}
The proof consists in the following computation.
    \begin{align*}
    0 &\stackrel{\eqref{idcolchy1y2}}{=}[y_1, z_1, z_2, y_2] \\
    &\stackrel{\eqref{[y,z,z]}}{=} y_1z_1z_2y_2+z_2z_1y_1y_2-y_2y_1z_1z_2-y_2z_2z_1y_1 \\
    & = y_1z_1z_2y_2+[z_2,z_1]y_1y_2+z_1z_2y_1y_2-y_1y_2z_1z_2-y_2[z_2,z_1]y_1 -y_2z_1z_2y_1 \\
    &\stackrel{\eqref{idcolchy1y2}, \eqref{2y2z}}{=}z_1z_2y_1y_2-y_1y_2z_1z_2. \qedhere 
    \end{align*}
\end{proof}

Inductively, for every integers $m$, $n \geq 0$ we have the following weak graded identity of the pair $(M, S)$.
\begin{equation}
   y_1y_2 \ldots y_{2m}z_1z_2 \ldots z_{2n} = z_1z_2 \ldots z_{2n}  y_1y_2 \ldots y_{2m} \label{y2mz2n=z2ny2m}.
\end{equation}

\section{
\texorpdfstring{A basis of weak graded identities for $(M, S)$ }{A basis of weak graded identities for (M, S)}}

In the previous section, we listed several weak $\mathbb{Z}_2$-graded identities of the pair $(M, S)$. In this section, we shall describe the form of the polynomials in $K\left<X \right>$ modulo $I$. We denote $F = K\left<X \right> / I$.

\begin{lemma} \label{lemma simetria ys}
    Let $f(y_1, y_2, \ldots , y_m, z_1, z_2, \ldots , z_n)$ be a multihomogeneous polynomial in $F$. Then $f$ is symmetric in its even variables $y_1$, $y_2$, \dots, $y_m$.
\end{lemma}
\begin{proof}
    It is sufficient to prove that, in an arbitrary monomial $m \in F$, we can rearrange the indices of the even variables as we want, keeping the coefficient (and the sign) of the monomial. Thus we will show that transposing two of the even variables in a monomial does not change it, modulo $I$.

Let $Y_i$ and  $Z_j$ denote words with only even and odd variables, respectively. If $m= Y_1Z_1$ or $m= Z_1Y_1$, where $Z_1$ is allowed to be empty, by the identity \eqref{idcolchy1y2}, we can reorder the indices of the even variables in the way we want.

Suppose now that $m$ has a factor of the form $y_i Z_k y_j$. We can assume that  $Z_k$ is of even length by \eqref{zimpar}. In this case, by \eqref{idyz2ky}, we can transpose $y_i$ and $y_j$.

    Now, let $m=Y_1Z_1Y_2Z_2 \ldots Y_kZ_k$ be a monomial, where some of the $Y_1$ or $Z_k$ may be empty. Combining the two observations above, we can rearrange the indices of the $y$'s in the way we want. This finishes the proof.
\end{proof}

Since the characteristic of the ground field is zero, this lemma implies that if $f (y_1, y_2, \ldots , y_m , z_1, z_2, \ldots , z_n) \in K\left<X \right>$ is a weak graded identity of the pair $(M, S)$ then $f$ is equivalent to a polynomial in a single even variable, $g (y, z_1, z_2, \ldots , z_n)$. In other words, we can symmetrize (restitute) the symmetric variables. Pay attention that the degree of $f$ in the new even variable $y$ is equal to the sum of the degrees of $f$ in $y_1$, $y_2$, \dots, $y_m$.

Now we study the symmetry of the polynomials $f \in F$ with respect to the odd variables $z \in Z$. Let $\tau_{ij} $ denote the transposition of the two odd variables $z_i$ and $z_j$ and let $\deg_x f $ denote the degree of $f$ with respect to the variable $x \in X$. 

\begin{proposition} \label{propf(y,z,z_1,z_2)}
    Let $f(y,z, z_1, z_2)$ be a multihomogeneous polynomial in $F$ which is linear in $z_1$ and $z_2$, $\deg  _z f>0$. Then $(1 - \tau _{12})f = \varphi (y, y_1, z)_{|y_1=[z_1,z_2]}$, where $\varphi$ is a suitable polynomial that is linear in the ``new'' variable $y_1$.
\end{proposition}
\begin{proof}
    It is sufficient to prove the statement for a monomial $m=z_1gz_2$, where $g$ is a word depending only on $y$ and $z$.

    First, we consider the monomial $g=y^{a_1}z^{b_1}y^{a_2}z^{b_2} \ldots y^{a_m}z^{b_m}$, $a_i$, $b_i \geq 0$. By the identity \eqref{zimpar}, if $a_i>0$ and $a_{i+1}>0$, then $b_i$ is even. From \eqref{yimpar}, we can assume that if $b_i >0$ and $b_{i+1} >0$ then $a_i $ is even. Since \eqref{y2z2=z2y2} is valid on $(M,S)$, after manipulations, we can suppose that $g$ has one of the following forms:
    \begin{equation} \label{polin em y e z}
      y^{r}, \hspace{0,4cm}z^s, \hspace{0,4cm}y^rz^s, \hspace{0,4cm}z^sy^r, \hspace{0,4cm}zy^{2r}z^{2s+1}, \hspace{0,4cm}y^rz^{2s}y, 
    \end{equation}
where $r$ and $s$ are non negative integers. 

If $g=y^{2k}$, by \eqref{zy2z}, $(1- \tau _{12})m =0$ modulo $I$.

If $g=y^{2k+1}$, we take into account the possible positions of $z$. By \eqref{idyzy}, in $m=g_1z_1y^{2k+1}z_2g_2$ we must have a variable $z$ immediately to the right of $z_2$ or immediately to the left of $z_1$. For $m= zz_1y^{2k+1}z_2 $, we have
\[
zz_1y^{2k+1}z_2=[z,z_1]y^{2k+1}z_2+z_1zy^{2k+1}z_2 \stackrel{\eqref{idcolchy1y2}, \eqref{yimpar}}{=}y^{2k+1}[z,z_1]z_2-z_1z_2y^{2k+1}z.
\]
In the last summand of the last equality, permuting $z_1$ and $z_2$ we achieve the statement, while the first summand is in the form where $g=z^s$. The case where $m=z_1y^{2k+1}z_2 z$ is dealt with analogously.

For $g=z^s$, we use induction on $s$.
If $s=1$, then $z_1zz_2 \stackrel{\eqref{s3embolinha}}{=}z \circ [z_1,z_2]$. For $s=2$,
\[
z_1z^2 z_2 =[z_1,z][z,z_2]+zz_1zz_2-zz_1z_2z+z_1zz_2z.
\]
Since $(1- \tau _{12})[z_1,z][z,z_2] \stackrel{\eqref{idcolchy1y2}}{=}0$, the case $s=2$ is reduced to the cases $s=0$ (trivial) and $s=1$ (treated above).

  For every $s>1$ we have that
  \begin{align*}
      (1-\tau_{12})z_1z^sz_2 = &(1-\tau_{12})([z_1,z]z^{s-2}[z,z_2])+ \\
      &+(1-\tau_{12})(zz_1z^{s-1}z_2+z_1z^{s-1}z_2z-zz_1z^{s-2}z_2z).
  \end{align*}
On the right-hand side, the first summand vanishes (a consequence of \eqref{zimpar} or \eqref{idyz2ky}) while in the remaining summands, the distance between $z_1$ and $z_2$ is less than $s$. By induction, we are done.

By applying \eqref{yimpar} if $r$ is odd, and \eqref{zy2z} if $r$ is even, we obtain that $z_1y^rz^sz_2 = \pm zy^rz_1z^{s-1}z_2$. Thus the case $g=y^rz^s$ is reduced to the case $g=z^s$. Similarly, we treat the case $g=z^sy^r$.

Now, let $g=zy^{2r}z^{2s+1}$. 
\begin{align*}
    z_1zy^{2r}z^{2s+1}z_2 &=z_1zy^{2r}z^{2s}[z,z_2]+z_1zy^{2r}z^{2s}z_2z \\ & \stackrel{\eqref{zy2z}}{=}z_1z[z,z_2]y^{2r-1}z^{2s}y+z_1zy^{2r}z^{2s}z_2z, 
\end{align*}
and we use induction and the previous cases.

If $g=y^rz^{2s}y$, then $m=z_1y^rz^{2s}yz_2=-z_1y^rz^{2s-1}z_2yz$ and we are in a position to apply a case discussed above. This completes the proof.
\end{proof}

\begin{remark} \label{f(y,z, z1, z2)=g(y,z)}
    Since $S$ is a simple Lie algebra, we have $[S, S] = S$. As $\mathrm{char} K = 0$, the last proposition implies that if $f(y,z, z_1, z_2)$ is a weak graded identity of $(M, S)$ that is linear in $z_1$ and $z_2$ and $\deg  _z f>0$, then the polynomial $\varphi$ of the statement is also a weak graded identity. 
    Moreover, by Lemma \ref{lemma simetria ys}, $f$ is a consequence of a polynomial depending only on one even variable and one odd variable. 
\end{remark}

We denote $\langle \mathrm{ad} a , \mathrm{ad} b , \mathrm{ad} c \rangle = s_3(\mathrm{ad} a , \mathrm{ad} b , \mathrm{ad} c)$. 
The identities in the next lemma can be proved by direct computation.

\begin{lemma}
    The following identities hold in every associative algebra
    \begin{equation}
        [a,b,c,d]-[a,b,d,c]-\left[[a,b], [c,d]\right]
    \end{equation}
    \begin{equation}
        \left< \mathrm{ad} a , \mathrm{ad} b , \mathrm{ad} c \right>d -  \left[ [d,a], [b,c]\right]- \left[ [d,b], [c,a]\right]- \left[ [d,c], [a,b]\right] \label{adcomocolcdecolc}
    \end{equation}
\end{lemma}

In order to prove a kind of symmetry for the odd variables, we need the following results. Some of them are based on results of \cite{trishin1983identities}.

\begin{lemma}
    Modulo the identities \eqref{idcolchy1y2}, \eqref{idyzy}, the following identity on the $\mathbb{Z}_2$-graded pair $(M,S)$ holds:
    \begin{equation}
        \left< \mathrm{ad} z_2 , \mathrm{ad} z_3 , \mathrm{ad} y \right>z_1 - 2 \left[ [z_1, y], [z_2,z_3]\right] \label{adzzy}.
    \end{equation}
 \end{lemma}

\begin{proof} 
A computation shows that
  \begin{align*}
         \left< \mathrm{ad} z_2 , \mathrm{ad} z_3 , \mathrm{ad} y \right>z_1 \stackrel{\eqref{adcomocolcdecolc}}{=} &[[z_1,z_2],[z_3,y]]+[[z_1,z_3],[y,z_2]]+[[z_1,y],z_2,z_3] \\
         = & [z_1,z_2]z_3y- [z_1,z_2]yz_3 -z_3y[z_1,z_2]+yz_3[z_1,z_2]+ [z_1,z_3]yz_2 \\
         & \, -[z_1,z_3] z_2y-yz_2[z_1,z_3] +z_2y[z_1,z_3]   +[[z_1,y],z_2,z_3] 
\\
          \stackrel{\eqref{idyzy}}{=} & -yz_1z_2z_3+yz_2z_1z_3-z_3z_1z_2y+z_3z_2z_1y+yz_1z_3z_2\\
         & \, -yz_3z_1z_2+z_2z_1z_3y-z_2z_3z_1y +[[z_1,y],z_2,z_3] \\
         = & yz_1[z_3,z_2]+yz_1[z_2,z_3]+y[z_2,z_3]z_1+z_1[z_2,z_3]y+ \\
         & \, +[z_2,z_3]z_1y  +[[z_1,y],z_2,z_3] \\
         \stackrel{\eqref{idyzy}}{=} &  2 \left[ [z_1, y], [z_2,z_3]\right]. \qedhere
  \end{align*}  
\end{proof}

\begin{lemma} \label{lemma3[]}
    Let $f = v_1z_iv_2z_jv_3$ be a polynomial where $z_i$ and $z_j$ are odd variables and among the commutators $v_1$,  $v_2$, $v_3$, we can find one $v_q=[u_1,u_2]$ of length greater than 1. Then, for some multilinear polynomial $\varphi$, the identity
\[
 (1-\tau_{ij} )f = \varphi (y, \hat{z}_i, \hat{z}_j, \ldots , \hat{v}_q, \ldots , u_1, u_2) _{|y=[z_i, z_j]}
\]
    is a consequence of the identities \eqref{idcolchy1y2}, \eqref{idyzy}, \eqref{idzyz}, and \eqref{s3}.
\end{lemma}
\begin{proof}
 Suppose $v_2 = [u_1,u_2]$. We will work with all the possible parities of $u_1$ and $u_2$.
 \begin{itemize}
     \item If both $u_1$, $u_2$ are odd, we have
     \begin{align*}
         (1-\tau_{ij} )f = &v_1(z_i [u_1,u_2]z_j - z_j[u_1,u_2]z_i)v_3 \\
        \stackrel{\eqref{s3emcolchetes}}{=}& v_1\{(u_1[z_i,u_2]-u_2[z_i,u_1])z_j-(u_1[z_j,u_2]u_2+[z_j,u_1])z_i\}v_3 \\
           \stackrel{\eqref{s3emcolchetes}}{=}& v_1(u_2[z_j, z_i]u_1+u_1[z_i,z_j]u_2)v_3.
     \end{align*}
     \item If both $u_1$, $u_2$ are even, then $f \equiv 0$ by \eqref{idcolchy1y2}.
     \item  If $u_1$ and $u_2$ have distinct parity, we can write, without loss of generality $u_1= z \in L \left< X \right> _1$ and $u_2 = y \in L\left< X \right> _0$.
     \begin{align*}
         (1-\tau_{ij}  )z_i[z,y]z_j = & z_izyz_j - z_iyzz_j- z_jzyz_i+ z_jyzz_i \\
         \stackrel{\eqref{idzyz}}{=}& -z_iz_jyz+zyz_iz_j+z_jz_iyz-zyz_jz_i \\
         = & zy[z_i,z_j] - [z_i,z_j]yz
     \end{align*}
      \end{itemize}
 Now suppose $v_1 = [u_1,u_2]$. Since $(1-\tau_{ij} )f = (1-\tau_{ij} ) v_1 z_iz_jv_3 =  v_1 [z_i,z_j]v_3$, we can suppose that $v_2 $ has length at least one.
 \begin{itemize}
     \item If $v_2$ is odd, by \eqref{s3embolinha}, we have
\[
(1-\tau_{ij} )f = (1-\tau_{ij} )[u_1,u_2]z_iv_2z_j = [u_1,u_2](v_2 \circ [z_j,z_i]).
\]
     \item If $v_1$ and $v_2$ are even, then $f=0$ by \eqref{idyzy}.
     \item If $v_1$ is odd and $v_2$ is even, we write $[u_1,u_2] = z$ and $v_2 = y$ in order to simplify the notation.
     \begin{align*}
        v_1z_iv_2z_j = &z z_i y z_j = z[z_i,y]z_j+zyz_iz_j \\
         = & [z,z_i]yz_j+z_izyz_j  +zyz_iz_j \\
         \stackrel{\eqref{idzyz}}{=} &[z,z_i,y]z_j+y[z,z_i]z_j+z_iz_jyz +zyz_iz_j\\
         \stackrel{\eqref{idcolchy1y2}}{=} & y[z,z_i]z_j+z_iz_jyz +zyz_iz_j
     \end{align*}
     and we can treat this case using the previous cases.
 \end{itemize}
 The case where $v_3=[u_1,u_2]$ is similar to the case where $v_1=[u_1,u_2]$.
\end{proof}

\begin{lemma} \label{lemmacomutadorlongo}
 Let $f = v_1 v_2 \cdots v_k z_iv_{k+1} \cdots v_{k+r}z_jv_{k+r+1} \cdots v_l$ be a multilinearpolynomial where $z_i$ and $z_j$ are odd variables and among the commutators $v_1$, $v_2$, \dots, $v_l$ we can find one $v_q=[u_1,u_2]$ whose length is greater than one. Then, for some multilinear polynomial $\varphi$, the identity
\[
 (1-\tau_{ij} )f = \varphi (y, \hat{z}_i, \hat{z}_j, \ldots , \hat{v}_q, \ldots , u_1, u_2) _{|y=[z_i, z_j]}
\]
is a consequence of the identities \eqref{idcolchy1y2}, \eqref{idyzy}, \eqref{idzyz}, and \eqref{s3}.
\end{lemma}
\begin{proof}
The idea of the proof consists of first reducing the distance between $z_i$ and $z_j$, and then applying induction and the previous Lemma.

\textbf{Case 1: } $q\geq k+r+1$ or $q \leq k$.

 By the equality $x_1x_2= [x_1,x_2]+x_2x_1$ and induction, we can suppose $q=k+r+1$ or $q=k$. First, we assume that $q=k+r+1$.
As above we consider the parities of the commutators. 
 
Note that, if $v_{k+r} $  and $v_q$ are both even, by the identity \eqref{idyzy},  $v_{k+r}z_jv_q=0$. If $v_{k+r} $ is even and $v_q$ is odd, by \eqref{[y,z,z]} and \eqref{[z_1, y, z_2] = [z_2, y, z_1]}, we have
\[
v_{k+r}z_jv_q =[v_{k+r}, z_j,v_q] -v_qz_jv_{k+r}=[v_{k+r}, v_q, z_j]-v_qz_jv_{k+r} .
\]
In this representation, $v_{k+r}z_jv_q$ is written as summands in the form that will be treated in Case 2, but the distance between $z_i$ and $z_j$ was reduced. 
 Similarly, if $v_{k+r}$ is odd and $v_q$ is even, 
  $v_{k+r}z_jv_q =[v_{k+r}, z_j,v_q] -v_qz_jv_{k+r}=[ v_q, v_{k+r}, z_j]-v_qz_j v_{k+r}$.
  
For the remaining case, we use 
\[
v_{k+r}z_jv_q = v_qz_jv_{k+r}-[v_q, v_{k+r}]z_j-z_j[v_q,v_{k+r}]+ s_3(v_{k+r}, z_j,v_q). 
\]
On the right-hand side of the last equality, in the second and third summands we reduced the distance between $z_i$ and $z_j$. The last summand vanishes due to \eqref{s3}.  The first summand will be treated in Case 2. 

The case $q=k$ is handled similarly.

    \textbf{Case 2:} $k<q<k+r+1$.
    
    Once again, by the equality  $x_1x_2= [x_1,x_2]+x_2x_1$ and induction, we have to deal only with the case $z_i g v_qv_{k+r}z_j$, where $g=v_{k+1} \cdots v_{k+r-2}$.

    If $v_q$ is odd and $v_{k+r}$ is even, by \eqref{idzyz}, $v_qv_{k+r}z_j=-z_jv_{k+r}v_q$ and we are back to  Case 1, but the distance between $z_i$ and $z_j$ decreased.

    Now we write 
\[
v_qv_{k+r}z_j = v_qz_jv_{k+r}+v_{k+r}z_jv_q+z_j[v_q,v_{k+r}]-z_jv_qv_{k+r}+[v_q, [v_{k+r}, z_j]].
\]
In this equality, we have reduced the distance between $z_i$ and $z_j$ in all summands, except the last one. We treat it separately, in accordance with the possible parities for $v_q$ and $v_{k+r}$.
    
    If $v_q$ is even and $v_{k+r}$ is odd, $[v_q, [v_{k+r}, z_j]]=0$, by  identity \eqref{idcolchy1y2}. 

    Suppose both $v_q $ and $v_{k+r}$ are odd. Then, since we assumed that $v_q$ is a nontrivial commutator, we can write $v_q = [y,z]$ for some $y \in L \left< X \right> _0 $ and some $z \in L \left< X \right> _1$. Using  the identity \eqref{adzzy} and the the fact that the element $\left< \mathrm{ad} x_1, \mathrm{ad} x_2, \mathrm{ad} x_3 \right>$, $x_1$, $x_2$, $x_3 \in X$, commutes with every element of the relatively free algebra $F$, we have 
\[
z_ig[v_q, [v_{k+r}, z_j]]= z_i g[[y,z],[v_{k+r}, z_j]]= (1/2)z_i\left< \mathrm{ad} v_{k+r} , \mathrm{ad} z_j , \mathrm{ad} y \right> gz, 
\]
 where the distance between $z_i$ and $z_j$ does not exceed 1.   

If $v_q$ and $v_{k+r}$ are even, we can assume $v_q = [z, \tilde{z}]$, where $z$,  $\tilde{z} \in L \left< X \right> _1$. Therefore
\[
z_ig[v_q, [v_{k+r}, z_j]]=z_ig[[z, \Tilde{z}],[v_{k+r},z_j]] \stackrel{\eqref{adzzy}}{=} (1/2) 
  z_i  \left<\mathrm{ad} v_{k+r} , \mathrm{ad} \Tilde{z} , \mathrm{ad} z_j  \right>gz.
\]
On the right-hand side of the last equality, the distance between $z_i$ and $z_j$ does not exceed 1.

   Using induction and the previous Lemma, it remains to prove the statement for the polynomial $[u, z_i, z_j, v]$. 

But $(1-\tau_{ij} )  [u, z_i, z_j, v] = [u, z_i, v]|_{z_i = [z_i, z_j]}$.   
\end{proof}

\begin{remark}
    We draw the readers' attention to the interplay between Cases 1 and 2. We need  Case 2 to work on  Case 1 and vice versa, but the distance between $z_i$ and $z_j$ does not increase in Case 1 while it decreases in Case 2. Therefore eventually we use the base of the induction.
\end{remark}

\begin{lemma}
    Let $f$ be a multilinear polynomial of degree $n \geq 4$ that is skew-symmetric in two different pairs of odd variables, and let us call them $z_i$, $z_j$ and $z_p$, $z_q$. Then there exists a multilinear polynomial $\varphi$ of degree $n-1$, such that the $\mathbb{Z}_2$-graded identity $f=\varphi _{|y=[z_i,z_j]}$ belongs to the $\mathbb{Z}_2$-graded weak ideal $I$.
\end{lemma}

\begin{proof}
    It suffices to prove the statement for the polynomial 
\[
f=(1-\tau_{ij} )(1-\tau _{pq})x_1 \cdots x_k z_i x_{k+1} \cdots x_{k+r}z_jx_{k+r+1} \cdots x_{n-2},
\]
where $z_p$, $z_q \in \{x_1, \ldots x_{n-2}\} \subseteq X = Y \cup Z$, and $\tau_{uv}$ denotes the transposition of the corresponding odd variables.

If $ k<p,q \leq k+r$, using $ab=[a,b]+ba$, we can suppose that 
\[
(1-\tau _{pq})x_1 \cdots x_k z_i x_{k+1} \cdots x_{k+r}z_jx_{k+r+1} \cdots x_{n-2}
\]
is a sum of polynomials containing nontrivial commutators. Then, by the previous Lemma, there exists  $\varphi$ as in the statement.
    The cases where $p$, $q \leq k$ or $p$, $q > k+r$ are analogous.

    If $p\leq k$ and $k<q \leq k+r$, once again by the equality $ab=[a,b]+ba$ and Lemma \ref{lemmacomutadorlongo}, it is sufficient to analyze the case where $p=k$ and $q=k+1$. But then, by \eqref{s3}, $(1-\tau _{pq})z_pz_iz_q =z_i[z_q,z_p]+[z_q,z_p]z_i$, and once more we apply  Lemma \ref{lemmacomutadorlongo} to conclude this part of the result.
    The cases where $k<p\leq k+r$ and $q>k+r$ are analogous.

    The remaining cases are similar to the previous ones, changing the role of the pairs $z_i$, $z_j$ and $z_p$, $z_q$.
\end{proof}

\begin{proposition} \label{prop simetria em n-2 zs}
    Let $f$ be a multilinear monomial of degree $k+n \geq 4$ depending on $k$ even variables and $n$ odd variables with $k\geq 0$, $n \geq 4$. Among the odd variables, we single out two variables, $z_i$ and $z_j$.
    Then, in the pair $(K \langle X \rangle , L \langle X \rangle)$, the polynomial $(1-\tau_{ij} )f$ can be written, modulo $I$, as
    \begin{equation}
        \frac{(1-\tau_{ij} )}{(n-2)!} \displaystyle \sum _{\sigma \in S_{n-2}} \sigma f + \varphi (y, y_1, \ldots , y_k, z_1, \ldots \hat{z_i}, \ldots , \hat{z_j}, \ldots , z_{n})_{|y=[z_i,z_j]} \label{simetria de n-2 zs}.
    \end{equation} 
Here $S_{n-2}$ are the permutations of the set $  \{ 1, \ldots , \hat{i}, \ldots , \hat{j}, \ldots , n\}$, $\varphi$ is a multilinear polynomial of degree $n+k-1$, and $\sigma f$ stands for
\[
f (y_1, \ldots , y_k ,  z_{\sigma (1)}, \ldots ,z_{\sigma (i-1)}, z_i, z_{\sigma (i+1)} ,\ldots ,z_{\sigma (j-1)}, z_j ,z_{\sigma (j+1)}, \ldots , z_{\sigma (n) }).
\]
\end{proposition}

   \begin{proof}
        Since every permutation may be represented as a product of transpositions, for every $\sigma \in S_{n-2}$, the polynomial $(1-\tau_{ij} )(1-\sigma)f$ satisfies the conditions of the previous lemma. Therefore, there exists a multilinear polynomial $\varphi _\sigma$ of degree $k+n-1$ such that $(1-\tau_{ij} )(1-\sigma)f = \varphi _{\sigma }{_{|y=[z_i,z_j]}}$. Then,
\[
 \sum _{\sigma \in S_{n-2}}(1-\tau_{ij} )(1-\sigma)f= \sum  _{\sigma \in S_{n-2}} \varphi _\sigma,
\]
and this in turn implies
\[
(n-2)!(1-\tau_{ij} )f-\sum _{\sigma \in S_{n-2}} (1-\tau_{ij} )\sigma f =  \sum  _{\sigma \in S_{n-2}} \varphi _\sigma.
\]
        Denoting $\varphi =  1 / (n-2)! \sum  _{\sigma \in S_{n-2}} \varphi _\sigma$, the proof is completed.
    \end{proof}

We recall that, since the base field is infinite, the multihomogeneous components of an identity are also identities. Then, as in Remark \ref{f(y,z, z1, z2)=g(y,z)}, we can conclude that the polynomial  $\varphi $ of \eqref{simetria de n-2 zs} is an identity of the $\mathbb{Z}_2$-graded pair $(M,S)$. 
 Inductively, 
     $\varphi$ is a linear combination of polynomials  which are symmetric in $l-2$  of its  $l \leq n-1$ odd variables. Eventually, we conclude that $\varphi$ is a consequence of polynomials of the form $\varphi _1(y, z_1, z_2)$ and $\varphi _2(y, z, z_1, z_2)$, and these are linear in $z_1$ and $z_2$. The latter of these polynomials, in turn, is a consequence of a weak graded identity of the form $g(y,z)$. 

     Likewise, the polynomial $(1-\tau_{ij} ) \sum _{\sigma \in S_{n-2}} \sigma f$ in \eqref{simetria de n-2 zs} is a weak graded identity too, which is symmetric with respect to $n-2$ odd variables and all its even variables. Therefore, it is equivalent to a weak graded identity of the form $g_1(y, z, z_1, z_2)$. Then, by Remark \ref{f(y,z, z1, z2)=g(y,z)}, $(1-\tau_{ij} ) \sum _{\sigma \in S_{n-2}} \sigma f$ is a consequence of a polynomial of the form ${g}(y,z)$.
     
Therefore we can suppose that $(1-\tau_{ij}  )f$ is a consequence of polynomials $g(y,z)$ and $\tilde{g}(y, z_1, z_2)$, where the latter polynomial is linear in $z_1$ and $z_2$. 
Since every $\sigma \in S_n$ can be written as a product of transpositions, we have that $f - \sigma f$ is a consequence of polynomials of the form  $g(y,z)$ and $\tilde{g}(y, z_1, z_2)$, for every $\sigma \in S_n$. Summing $f - \sigma f$ over all $\sigma \in S_n$ we obtain the following corollary.

\begin{corollary} \label{cor simetria dos z's}
   Let $f(y_1, y_2, \ldots , y_m , z_1, z_2, \ldots , z_n)=0$ be a multilinear  $\mathbb{Z}_2$-graded identity of the pair $(M,S)$, with $n\geq 3$. Then $f$ is a consequence of  polynomials of the form ${g}(y,z)$ and $\tilde{g}(y, z_1, z_2)$, linear in $z_1$ and $z_2$.  
\end{corollary}

We will need the following analogues of the well known generic matrix algebras.
Let $ A := K [a_i, b_i, c_i \mid i \geq 1]$ be the commutative polynomial algebra generated by the independent variables $a_i$, $b_i$, $c_i$, $i \geq 1$. Consider the associative algebra $R$ and the Lie algebra $L$ generated by the generic matrices 
\begin{displaymath}
Y_i=\left(\begin{array}{ccc}
0 & 0 & 0 \\
0 & 2a_i & 0 \\
0 & 0 & -2a_i
\end{array}\right) \, \, \, \, \mathrm{and}  \, \, \, \, Z_i= \left(\begin{array}{ccc}
0 & -b_i & c_i \\
-2c_i & 0 & 0 \\
2b_i & 0 & 0
\end{array}\right), \hspace{0,4cm} i\geq 1.
\end{displaymath}
Observe that $R$ and $R^{(-)}$ are $\mathbb{Z}_2$-graded associative and Lie algebras, respectively, with the $\mathbb{Z}_2$-grading defined at the very beginning of the previous section. Also, $L$ is a $\mathbb{Z}_2$-graded Lie algebra, with $L_i= R^{(-)}_i \cap L$, $i=0$, 1.

Then a standard argument shows that the $ \mathbb{Z}_2$-graded pair $(R, L)$ is the relatively free $\mathbb{Z}_2$-graded pair of countable rank in
the variety of $\mathbb{Z}_2$-graded pairs generated by $T_2(M, S)$  (see for example \cite[Lemma 3]{koshlukov2012basis}). Therefore we shall work in the $\mathbb{Z}_2$-graded pair $(R, L)$ when necessary.

Recall that the weak $T_2$-ideal $I$ is contained in $T_2 (M, S)$ and, as stated in the last corollary, the elements in $F= K\left<X \right> /I$ are consequences of polynomials in one even variable $y \in Y$ and at most two odd variables $z_1$, $z_2 \in Z$. We draw the readers' attention to the fact that these polynomials are multihomogeneous but they need not be linear, apart from the case of two odd variables where the polynomial is linear in both $z_1$ and $z_2$.

\begin{proposition}
There are no nontrivial weak graded identities for the pair $(M, S)$ modulo the $T_2$-ideal $I$. 
\end{proposition}

\begin{proof}
  By Corollary \ref{cor simetria dos z's}, if we want to study the weak graded identities for $(M, S)$ modulo $I$ it is sufficient to analyze the polynomials of the form $g(y,z)$ and $g(y, z_1, z_2)$, the latter is linear in $z_1$ and $z_2$, since there are no weak graded identities depending only on $y$ (that is of the form $y^m$) or only on $z$ (that is of the form $z^n$), modulo $I$.
  
Thus, in order to prove that there is no nontrivial weak graded identity for the pair $(M, S)$ modulo $I$, it is enough to prove that the monomials in $K\left<X \right> /I$ depending on $y$ and $z$ or on $y$, $z_1$ and $z_2$ are linearly independent. To achieve this, in turn, we prove the linear independence for every set of multihomogeneous monomials of the same multidegree.

We work in the $\mathbb{Z}_2$-graded pair $(R,L)$, setting $y = Y$ and $z_i = Z_i$, where $Y$ and $Z_i$ are the matrices introduced above. 

Let us analyze first, according to their multihomogeneous degree, the polynomials modulo $I$ that depend only on $y$ and $z$, which are, as we have seen above, $ y^rz^s$, $z^sy^r$, $zy^{2r}z^{s-1}$, and $ y^rz^{2s}y$.

\begin{itemize}
    \item If $\deg _y f =r$ and $ \deg_z f =2k$ are both even, the possibilities for a monomial with such multidegree are $y^rz^{2k}, \, zy^rz^{2k-1}$ and $y^{r-1}z^{2k}y$. Then evaluating on the generic matrices we have 
    \begin{displaymath}
     Y^rZ^{2k}=2^{r+2k-1}\left(\begin{array}{ccc}
0 & 0 & 0 \\
0 & a^rb^kc^k & -a^rb^{k-1}c^{k+1} \\
0 & -a^rb^{k+1}c^{k-1} & a^rb^kc^k
\end{array}\right),
    \end{displaymath}
    \begin{displaymath}
        Z Y^{r}Z^{2k-1}=2^{r+2k+2}\left(\begin{array}{ccc}
a^rb^kc^k & 0 & 0 \\
0 & 0 & 0 \\
0 & 0 & 0
\end{array}\right),
    \end{displaymath}
    \begin{displaymath}
        Y^{r-1}Z^{2k}Y=2^{r+2k}\left(\begin{array}{ccc}
0 & 0 & 0 \\
0 & a^rb^kc^k & a^rb^{k-1}c^{k+1} \\
0 & -a^rb^{k+1}c^{k-1} & a^rb^kc^k
\end{array}\right).
    \end{displaymath}
By the form of the above matrices, it is clear that the second matrix cannot participate in a linear combination. Then a linear combination of the first and the third matrix equals 0 if and only if it is the trivial combination. Therefore no nontrivial linear combination of these elements is an identity.
    
    \item If $\deg _yf  =r$ and $ \deg_z f =2k+1$ are odd, the possible monomials are $y^rz^s$ and $z^sy^r$. But
    \begin{displaymath}
      Y^rZ^s = 2^{r+2k+1}  \left(\begin{array}{ccc}
0 & 0 & 0 \\
-a^rb^{k}c^{k+1} & 0 & 0 \\
-a^rb^{k+1}c^k & 0 & 0
\end{array}\right),
    \end{displaymath}
    and
      \begin{displaymath}
      Z^sY^r = 2^{r+2k}  \left(\begin{array}{ccc}
0 & -a^rb^{k+1}c^k & -a^rb^{k}c^{k+1} \\
0 & 0 & 0 \\
0 & 0 & 0
\end{array}\right).
    \end{displaymath}
      Clearly, no nontrivial linear combination of these elements is a weak graded identity.
    \item If $\deg _y  f=r$ is odd and $ \deg_z f =2k$ is even, we have the following possibilities
     \begin{displaymath}
     Y^rZ^{2k} = 2^{r+2k-1}  \left(\begin{array}{ccc}
0 &0 &0 \\
0 & a^rb^kc^k & -a^rb^{k-1}c^{k+1} \\
0 & a^rb^{k+1}c^{k-1} & -a^rb^kc^k
\end{array}\right)
    \end{displaymath}
    and
        \begin{displaymath}
     Z^{2k}Y^r = 2^{r+2k}  \left(\begin{array}{ccc}
0 &0 &0 \\
0 & a^rb^kc^k & a^rb^{k-1}c^{k+1} \\
0 & -a^rb^{k+1}c^{k-1} & -a^rb^kc^k
\end{array}\right).
    \end{displaymath}
      No nontrivial linear combination of these elements is an identity.
\item If $\deg_y  f =r$ is even and $ \deg_z f =2k+1$ is odd, the  possibilities for monomials with such multidegrees are 
        \begin{displaymath}
     Y^rZ^{2k+1} = 2^{r+2k+1}  \left(\begin{array}{ccc}
0 &0 &0 \\
-a^rb^{k}c^{k+1} &0 &0 \\
a^rb^{k+1}c^{k} &0 &0 \\
\end{array}\right),
    \end{displaymath}
       \begin{displaymath}
     Z^{2k+1}Y^r = 2^{r+2k+1}  \left(\begin{array}{ccc}
0 & -a^rb^{k+1} c^k & a^rb^{k}c^{k+1} \\
0 &0 &0 \\
0 &0 &0 \\
\end{array}\right).
    \end{displaymath}
Clearly, there is no nontrivial linear combination in this case either.
\end{itemize}

Now we deal with the polynomials $f(y,z_1,z_2)$, which are linear in $ z_1$ and $z_2$. First, we examine the possible monomials $m$ that can appear in these polynomials.
From \eqref{idyzy}, there is no factor of the type $y z_iy$ in $m$. 
This fact, along with identity \eqref{y2mz2n=z2ny2m} implies   that $m$ can only be in one of the following forms:
$ z_iz_jy^{2r}$, $y^{2r-1}z_iz_jy$, and $ z_i y^{2r}z_j$ if the total degree of $f$ is even, or
$y^{2r+1}z_iz_j$, $z_iz_jy^{2r+1}$, and $z_iy^{2r+1}z_j$ otherwise, where $ i \neq j \in \{1,2\}$. 

Following the same argument as above, we have:
\begin{itemize}
    \item If $\deg f = 2r+2$,
    \begin{displaymath}
Z_iZ_jY^{2r} =  2^{2r+1}   \left(\begin{array}{ccc}
0 &0 &0 \\
0 & a^{2r}b_jc_i & -a^{2r}c_ic_j \\
0 & -a^{2r}b_ib_j & a^{2r}b_ic_j
\end{array}\right) ,
\end{displaymath}
\begin{displaymath}
Y^{2r-1}Z_iZ_jY =  2^{2r+1}   \left(\begin{array}{ccc}
0 &0 &0 \\
0 & a^{2r}b_jc_i & a^{2r}c_ic_j \\
0 & a^{2r}b_ib_j & a^{2r}b_ic_j
\end{array}\right) , 
\end{displaymath}
\begin{displaymath}
Z_iY^{2r}Z_j =  2^{2r+1}   \left(\begin{array}{ccc}
a^{2r}(b_ic_j+b_jc_i) &0 &0 \\
0 & 0 & 0 \\
0 & 0 & 0
\end{array}\right). 
\end{displaymath}
Since $i \neq j \in \{1,2\}$, we have the following linear combinations resulting zero: $Y^{2r}Z_iZ_j - Y^{2r}Z_jZ_i-Y^{2r-1}Z_iZ_jY+Y^{2r-1}Z_jZ_iY=0$ and $Z_iY^{2r}Z_j - Z_jY^{2r}Z_i = 0$, but these are consequences of \eqref{idcolchy1y2} and \eqref{zy2z}, respectively. Consequently, they are trivial in $R$ modulo $I$.

\item If $\deg f = 2r+1$,
\begin{displaymath}
Y^{2r+1}Z_iZ_j =  2^{2r+2}   \left(\begin{array}{ccc}
0 &0 &0 \\
0 & a^{2r+1}b_jc_i & -a^{2r+1}c_ic_j \\
0 & a^{2r+1}b_ib_j & -a^{2r+1}b_ic_j
\end{array}\right) ,
\end{displaymath}
\begin{displaymath}
Z_iZ_jY^{2r+1} =  2^{2r+2}   \left(\begin{array}{ccc}
0 &0 &0 \\
0 & a^{2r+1}b_jc_i & a^{2r+1}c_ic_j \\
0 & -a^{2r+1}b_ib_j & -a^{2r+1}b_ic_j
\end{array}\right), 
\end{displaymath}
\begin{displaymath}
Z_iY^{2r+1}Z_j =  2^{2r+2}   \left(\begin{array}{ccc}
 a^{2r+1}(b_ic_j-b_jc_i)&0 &0 \\
0 & 0 & 0 \\
0 & 0 & 0
\end{array}\right).
\end{displaymath}
Here we have $Y^{2r+1}Z_iZ_j - Y^{2r+1}Z_jZ_i-Z_iZ_j Y^{2r+1}-Z_jZ_iY^{2r+1} =0$ and $Z_iY^{2r+1}Z_j+Z_jY^{2r+1}Z_i=0$, which are trivial in $R$ modulo $I$ by identities \eqref{idcolchy1y2} and \eqref{yimpar}.
\end{itemize}
\end{proof}

Therefore we obtain the proof of the main result of this paper, summarized in the following theorem.

\begin{theorem}
    Let $K$ be a field of characteristic zero and let $(M,S)$ be the $\mathbb{Z}_2$-graded pair relative to the adjoint representation of $\mathfrak{sl}_2 (K)$ equipped with the only nontrivial $\mathbb{Z}_2$-grading on $\mathfrak{sl_2}(K)$, and the induced grading on  $M_3(K)$. Then the following graded identities form a basis for the $\mathbb{Z}_2$-graded identities of the pair $(M, S)$.
$$ [y_1, y_2], \hspace{0,4cm}
        y_1zy_2, \hspace{0,4cm}
        z_1 y z_2 + z_2 y z_1 , \hspace{0,4cm}
    s_3(z_1, z_2, z_3), $$
$$
  y_1 z_1 z_2 y_2 - y_2 z_1 z_2 y_1, \hspace{0,4cm}
        z_1y_1y_2z_2-z_2y_1y_2z_2 ,\hspace{0,4cm}
     s_3 (z_1,z_2, y) x - x s_3  (z_1,z_2, y) ,
    $$
where $x = y_i$, $z_i$.
\end{theorem}

One can naturally ask the following questions.

\begin{enumerate}
    \item 
    Let $K$ be algebraically closed and of characteristic 0, and let $\rho\colon \mathfrak{sl}_2(K)\to \mathfrak{gl}(V)$ be an irreducible representation of dimension $m+1$. The natural $\mathbb{Z}_2$-grading on $\mathfrak{sl}_2(K)$ induces a $\mathbb{Z}_2$-grading on $V$ as follows. Take a maximal vector $v_0\in V$ and form the canonical basis $v_0$, $v_1$, \dots, $v_m$ of $V$. Then $V$ is $\mathbb{Z}_2$-graded by assigning degree 0 to the $v_{2i}$ and degree 1 to the $v_{2i+1}$. This induces a $\mathbb{Z}_2$-grading on $M_{m+1}(K)$. 

    Describe a basis of the $\mathbb{Z}_2$-graded identities for $\rho$. 
    \item 
    If $K$ is of positive characteristic $p$, the irreducible representations of $\mathfrak{sl}_2(K)$ are finite dimensional. 

    Describe their $\mathbb{Z}_2$-graded identities.
\end{enumerate}

\end{document}